\documentclass[12pt,a4paper]{article}

\usepackage{amsmath,amssymb,latexsym,amsthm}
\usepackage[english]{babel}
\usepackage{eufrak}

\begin{document}

\newtheorem*{Dfn}{Definition}
\newtheorem{Theo}{Theorem}
\newtheorem{Lemma}[Theo]{Lemma}
\newtheorem{Prop}[Theo]{Proposition}
\newtheorem{Cor}[Theo]{Corollary}
\newcommand{\Pro}{\noindent{\em Proof. }}
\newcommand{\Rem}{\noindent{\em Remark. }}

\medskip
\noindent

\title{Loops which are semidirect products of groups}
\author{\'A. Figula (Erlangen,Debrecen) and K. Strambach (Erlangen)}
\date{}
\maketitle
\footnotetext[1]{This paper was supported by DAAD.}
\footnotetext[2]{{\em Key words and phrases:} splitting loops extensions of 
groups, sharply transitive sections, loops in groups of affinities, 
eigenvalues of products of matrices, Akivis algebras.}
\footnotetext[3]{2000 {\em Mathematics Subject Classification:} 20N05, 51A99,
 22A30, 15A18, 17D99.}
\small{Dedicated to P\'eter T. Nagy on the occasion of his 
60th birthday, in friendship.}

\begin{abstract} 
We construct loops which are semidirect products of groups of affinities.  
As their elements in many cases one may take  transversal subspaces of an 
affine space. In 
particular we obtain in this manner smooth loops having  Lie groups of 
affine real transformations as the groups generated by  left 
translations, 
whereas the groups generated by  right translations are  smooth groups of 
infinite dimension. We also determine the Akivis algebras of these loops.  
\end{abstract}

\bigskip
\centerline{\bf 0. Introduction}

\bigskip

In \cite{birken2}, \cite{birken1} and \cite{winterroth} constructions of 
proper loops are discussed which are
semidirect products of groups. Whereas 
in \cite{birken2} there are few  constructions of such loops related to 
Example 3 in 
\cite{bruck}, p. 128, in \cite{birken1} a general theory for 
loops which are 
semidirect products of groups is developed. In \cite{winterroth} examples of 
proper analytic Bol loops are presented  which are \glqq twisted semidirect 
products\grqq \  of two Lie groups.  

In our paper we show that a wide class of proper loops $L$ can be 
represented 
within the group of affinities of an affine space ${\mathcal A}$ of 
dimension $2n$ over a commutative field $\mathbb K$. They are semidirect 
products of groups of translations of  ${\mathcal A}$ by suitable subgroups 
$\Gamma _0'$   
of $GL(2n, \mathbb K)$. For many  of them we may take as elements  affine 
$n$-dimensional transversal subspaces of  ${\mathcal A}$. This 
representation of the loops $L$ depends in an essential manner on the 
existence of a regular orbit in the hyperplane at infinity of  
${\mathcal A}$ for the group $\Gamma _0'$. 

To realize our examples it is important to know the 
eigenvalues for certain  products of matrices in 
$GL(n, \mathbb K)$. Since there is no unique procedure for the calculation 
of the eigenvalues of the product  $AB$ from the eigenvalues of the matrices 
$A$ and $B$ we have devoted Section 3 to this problem and give  answers 
in special cases. 

If the field $\mathbb K$ is a topological field then we obtain topological 
loops, for  real or complex numbers the constructed loops are analytic. 
For smooth proper loops obtained in this paper the group  
topologically generated by the left translations is a Lie group.
The difference between the multiplication for semidirect products in our 
paper and the multiplication for \glqq twisted semidirect products\grqq \ 
 in \cite{winterroth} seems to be negligible. But the groups 
topologically generated by all  translations of analytic loops treated in 
\cite{winterroth} are Lie groups, whereas the analytic loops considered here have smooth 
transformation groups of infinite dimension as the groups 
generated by all translations. Moreover, we prove that already the groups  
topologically generated by the right   translations of these  loops are 
smooth groups 
having a normal abelian subgroup of infinite dimension. 

The Akivis algebras of smooth loops contructed in this paper are 
 semidirect products of Lie algebras. Moreover, there are non-connected  
proper smooth  loops having Lie algebras as their Akivis algebras.

\bigskip
\noindent

\centerline{\bf 1. Some basic notions of loop theory} 

\bigskip
\noindent
A set $L$ with a binary operation $(x,y) \mapsto x \cdot y$ is called a loop 
if there exists an element $e \in L$ such that $x=e \cdot x=x \cdot e$ holds 
for all $x \in L$ and the equations $a \cdot y=b$ and $x \cdot a=b$ have 
precisely one solution which we denote by $y=a \backslash b$ and $x=b/a$. 
The left 
translation $\lambda _a: y \mapsto a \cdot y :L \to L$ as well as the right translation 
$\varrho _a:y \mapsto y \cdot a:L \to L$ 
is a bijection of $L$ 
for any $a \in L$. 
\newline
A loop $(L, \cdot )$ is a semidirect product of $H$ by $K$ if $H$ and $K$ 
are 
subloops of $(L, \cdot )$ such that: (i) $H$ is a normal subloop of  
$(L, \cdot )$,  (ii) $L=H K$,  (iii) $H \cap K=\{ e \}$, where $e$ is the 
identity of  $(L, \cdot )$ (cf. \cite{birken1}, p. 81).   

Let $L$ be a topological space,  respectively a  
$\mathcal C^{\infty }$-differentiable manifold. Then 
 $(L, \cdot )$ is a  topological, respectively a differentiable loop if  the 
maps  
$(x,y) \mapsto x \cdot y$, 
$(x,y) \mapsto x \backslash y$, 
$(x,y) \mapsto y / x: L^2 \to L$ are continuous, respectively differentiable.

To any differentiable loop $L$ we may associate an Akivis algebra 
which is realized in the tangent space of $L$ at the identity $e \in L$ and which plays a 
  similar role as the Lie algebra in the case of a Lie group (cf. \cite{akivis}, 
\cite{hofmann}). 
An Akivis algebra $(A,[.,.], \langle .,.,. \rangle)$ is a real vector space with a bilinear 
skew-symmetric map $(x,y) \mapsto [x,y]: A \times A \to A$ (called the commutator map) and a 
trilinear map  $(x,y,z) \mapsto \langle x,y,z \rangle: A \times A \times A \to A$ (called the associator 
map) such that the following identity (called the Akivis identity) holds: 
\[ \big [ [x,y],z \big]+ \big [ [z,x],y \big]+ \big [ [y,z],x \big] = \] 
\[ \langle x,y,z \rangle +
 \langle y,z,x \rangle +  \langle z,x,y \rangle - \langle z,y,x \rangle - \langle x,z,y \rangle -  \langle y,x,z \rangle \] 
for all $x,y,z \in A$. 

Let $(A_1,[.,.], \langle .,.,. \rangle)$ and $(A_2,[.,.], \langle .,.,. \rangle)$ be 
Akivis algebras. A homomorphism $\alpha :A_1 \to A_2$ is a linear map such that 
$[x,y]^{\alpha }=[x^{\alpha }, y^{\alpha }]$ and 
$ \langle x,y,z \rangle ^{\alpha }= \langle x^{\alpha },y^{\alpha },z^{\alpha } \rangle$ for 
all $x,y,z \in A_1$ holds. 

The Akivis algebra $A$ is a semidirect product $A=N \rtimes M$ of a Lie 
algebra $N$ by a Lie 
algebra $M$ if there exist in $A$ Lie subalgebras $N$ and $M$ together with 
an  endomorphism 
$\alpha :A \to M$ such that $N$ is the kernel of $\alpha $ and the vector space $A$  is 
the direct sum $A=N \oplus M$.

Let $G$ be the group generated by the left 
translations of $L$ and let $H$ be the stabilizer of $e \in L$ in the group 
$G$. 
The left translations of $L$ form a subset of $G$ acting on the cosets 
$x H, x \in G$, such that for any given cosets $aH$ and $bH$ there exists 
precisely one left translation $ \lambda _z$ with $ \lambda _z a H=b H$. 
\newline
Conversely, let $G$ be a group, let H be a subgroup of $G$ and let $\sigma : G/H \to G$ be a 
section with $\sigma (H)=1 \in G$ such that the subset $\sigma (G/H)$ 
generates $G$ and 
acts sharply transitively on the factor space $G/H=\{x H; x \in G\}$ 
(cf. \cite{loops}, p. 18).  We call such a section sharply 
transitive. 
Then the multiplication  defined by 
$x H \ast y H=\sigma (x H) y H$ on  $G/H$ or by $x \ast y=\sigma(xyH)$ 
on $\sigma (G/H)$ yields a loop $L(\sigma )$. If $N$ is the largest normal 
subgroup of 
$G$ contained in $H$ then the factor group $G/N$ is isomorphic to the group 
generated by 
the left translations of $L(\sigma )$. 

The loop $L$ is a group if and only if the set $\{ \lambda_x; \ x \in L \}$ 
of left 
translations  is a group. A loop $L$ has the left inverse, 
respectively the right  inverse property 
if the identity $x^{-1} (x y)=y $, respectively $(y x) x^{-1}=y$ holds  for 
all $x,y \in L$.

\bigskip
\noindent
\centerline{\bf 2.  A general construction}

\bigskip
\noindent
Let $\Gamma _0$ be a subgroup of the general linear group $GL(n, \mathbb K)$ 
which is different from the identity $I$ and acts on the $n$-dimensional  
vector space 
$\mathbb K^n$ of column vectors, where 
$\mathbb K$ is a commutative field. 
Let $\Gamma $ be the group of matrices 
\[  g(u,v,M)= \left ( \begin{array}{ccc}
1 & 0 & 0 \\
u & M & 0 \\
v & 0 & M^{\delta } \end{array} \right ),  u,v \in  \mathbb K^n,\ M \in \Gamma _0,  \]
where $\delta :\Gamma _0 \to \Gamma _0$ is an endomorphism. 

\begin{Theo} Let ${\Xi }=T_B \Gamma _0'$ be the complex product of the 
groups 
\[ T_B=\{g(u, B u, I);\ u \in \mathbb K^n \} \ \hbox{with} \ B \in GL(n, \mathbb K) \ 
\hbox{and} \  \Gamma _0'=\{ g(0,0,M);\ M \in \Gamma _0 \}. \]
 If no element of the set $\{ 
M^{-1}B^{-1} M^{\delta } A;\ M \in \Gamma_0 \}$ has eigenvalue $1$ then 
for the subgroup  
$H=\{g(u,A u,I);\ u \in \mathbb K^n \} \in \Gamma$ with 
$A \in GL(n, \mathbb K)$  the mapping  $\sigma :\Gamma /H \to \Gamma $ 
defined by 
\[ g(0,v,M) H \mapsto g(M (B M-M^{\delta } A)^{-1} v, B M (B M-M^{\delta } A)^{-1} v,  M) \]
 is a sharply transitive section with image $\Xi $.  

The loop $L_{\Xi }$ 
corresponding to  $\sigma $ is isomorphic to the semidirect product 
$\mathbb K^n \rtimes \Gamma _0$ of 
 the normal subgroup $\mathbb K^n$ by the group  $\Gamma _0$ under the 
multiplication defined by
\[ (u_1,M_1) \ast (u_2,M_2)=(u_1+u_2^{\psi }, M_1 M_2)  \]
for all $u_i \in \mathbb K^n$, $M_i \in \Gamma _0$, where $\psi $ is the 
invertible linear map 

\medskip
\noindent
(I) \centerline{$x \mapsto [B-(M_1M_2)^{\delta }A (M_1M_2)^{-1}]^{-1}[M_1^{\delta }(B-M_2^{\delta }A M_2^{-1})] x$. }

\medskip
\noindent
The loop $L_{\Xi }$ is a group if and only if  the following property 

\medskip
\noindent
(II) \centerline{ $\{ 
M^{-1}B^{-1} M^{\delta } ;\ M \in \Gamma_0 \}=\{B^{-1}  \} $ } 

\medskip
\noindent
holds. This is equivalent to the condition that  
$T_B$ is a normal subgroup of $\Gamma $.

Let $\Delta '$ be the semidirect product 
of the group $T'_B$, which is the minimal normal subgroup of 
 $ \Gamma $ containing $T_B$, by the group $ \Gamma _0 ' $. 
If  in  $ \Gamma _0 $ there is an element $M_0$ such that the matrix 
$M_0^{-1}B^{-1} M_0^{\delta } B$ has no eigenvalue $1$ then one has 
$T'_B=\{ g(u,v,I);\ u,v \in \mathbb K^n \}$ and $\Delta '=\Gamma $.

Let $ \Theta $ be  the maximal normal subgroup of $\Gamma $ 
contained in $H$. Then $\Theta $ consists of matrices 
$g(v,Av,I),\ v \in V$, 
where $V$ is the maximal subspace of $\mathbb K^n$ for which 
$M^{\delta } v=(A M A^{-1}) v$ for all $v \in V$ and $M \in \Gamma _0$ 
holds.  The group 
$\Delta $ 
 generated by the set $\Xi $ of left translations  is the factor group 
$\Delta '/ (\Theta \cap \Delta ')$.  
  The stabilizer of the identity of  $L_{\Xi }$ is the group 
$ H / (\Theta  \cap \Delta ')$.

No loop $L_{\Xi }$ satisfies the left inverse as well as 
the right inverse property.

If $\mathbb K$ is a topological field then any loop 
 $L_{\Xi }$ is a topological loop.  For real or 
complex  numbers  any loop $L_{\Xi }$ is  analytic.

\end{Theo}
\Pro  Any matrix $g(x,y,M)$ has a unique decomposition 
\[ g(x,y,M)=g(0,y-M^{\delta }A M^{-1}x, M) g(M^{-1}x, A M^{-1}x,I). \]
Therefore the set $\{g(0,v,M);\ v \in \mathbb K^n, M \in \Gamma _0 \}$ 
forms a system of representatives of the left cosets of $H$ in $\Gamma $. 
A  mapping $\sigma : \Gamma /H \to \Xi $  is a sharply transitive section if and only if  
for given $g(0,v_1,M_1)$ and   $g(0,v_2,M_2)$  there is precisely one matrix 
$g(u, B u, M)$ in $\sigma(\Gamma /H)$ and one matrix $g(z, A z,I)$ in $H$ 
such that 

\medskip
\noindent
$(1)$ \centerline{$g(u, B u, M)g(0,v_1,M_1)=g(0,v_2,M_2)g(z, A z,I)$}

\medskip
\noindent
holds. Since $M^{-1}_2B^{-1} M^{\delta }_2 A $ has no eigenvalue $1$ the 
matrix $B M_2( I- M^{-1}_2B^{-1} M^{\delta }_2 A)$ $=$ $B M_2- M^{\delta }_2 A$
is invertible and the unique solution of  $(1)$ is given by 
\[ M=M_2 M_1^{-1},\ z=(B M_2- M^{\delta }_2 A)^{-1}(v_2-M_2^{\delta } M_1^{-\delta } v_1), \]
\[ u=M_2(B M_2- M^{\delta }_2 A)^{-1}(v_2-M_2^{\delta } M_1^{-\delta } v_1). \] 
The set $\Xi =T_B \Gamma _0'$ is the set of the left translations of the loop $L_{\Xi }$ 
defined by the sharply transitive section $\sigma $.

For the elements $g(u_i,B u_i,M_i) \in \sigma(G/H)$, $i=1,2$, one has 
\[ g(u_1,B u_1, M_1)g(u_2,B u_2, M_2) \in g(z,B z, M') H  \] 
with a unique $z \in \mathbb K^n$ and a unique  $M' \in \Gamma _0$.  
This yields  $M'=M_1 M_2$ and 
\[ [M_1^{\delta }(B-M_2^{\delta }A M_2^{-1})] u_2= [B-(M_1M_2)^{\delta }A (M_1M_2)^{-1}] (z-u_1). \]  The matrix 
$B-(M_1M_2)^{\delta }A (M_1M_2)^{-1}$ is invertible since 
\[ B-(M_1M_2)^{\delta }A (M_1M_2)^{-1}=
B[(M_1M_2-B^{-1} (M_1M_2)^{\delta }A) (M_1M_2)^{-1}]= \] 
\[ B[(M_1M_2(I-(M_1M_2)^{-1} 
B^{-1} (M_1M_2)^{\delta }A)) (M_1M_2)^{-1}] \]
and  $(M_1 M_2)^{-1} B^{-1} (M_1 M_2)^{\delta }A$ has no eigenvalue $1$. 
It follows  

\bigskip
\noindent
$(2)$ \centerline{$z=u_1+ [B-(M_1M_2)^{\delta }A (M_1M_2)^{-1}]^{-1}[M_1^{\delta }
(B-M_2^{\delta }A M_2^{-1})] u_2$. }

\bigskip
\noindent
Hence (cf. \cite{loops}, p. 18)  
the loop $L_{\Xi }$ is isomorphic to the loop defined on $\Xi $ by the 
multiplication 
\[g(u_1,B u_1, M_1) \circ  g(u_2,B u_2, M_2)=g(z,B z, M_1 M_2), \] 
where $z$ is given in $(2)$. Moreover,  $L_{\Xi }$ is 
isomorphic 
to the loop $\tilde{L}_{\Xi }$ defined on ${\mathbb K}^n \Gamma _0$ by 
$(u_1,M_1) \ast (u_2,M_2)=(u_1+  u_2^{\psi }, M_1 M_2)$, 
where   $\psi $ is a linear map given in $(I)$. 
Since $({\mathbb K}^n,I)$ is a normal subgroup of $\tilde{L}_{\Xi }$ the 
loop 
$\tilde{L}_{\Xi }$ is a semidirect product of the group ${\mathbb K}^n$ by 
the group $\Gamma_0$.

The loop $L_{\Xi }$  is a group if and only if the set $\Xi $ is a group. 
This is equivalent to the fact that $T_B \Gamma _0'= \Gamma _0' T_B $ 
or that  for given $g(u,Bu,I)$ and $ g(0,0,M)$ there are elements $g(u',Bu',I)$ and $g(0,0,M')$ such that 
\[ g(u,Bu,I)  g(0,0,M)=g(0,0,M') g(u',Bu',I). \] 
This is the case if and only if $M=M'$ and $B=M^{\delta }BM^{-1}$ which is 
equivalent to  $\{ 
M^{-1}B^{-1} M^{\delta } ;\ M \in \Gamma_0 \}=\{B^{-1}  \} $ or to the normality of $T_B$ 
in $\Gamma $.

If the  set $\Xi $ is not a group then  
$\Xi =T_B \Gamma _0' \neq  \Gamma _0' T_B$ and 
the  loop  $L_{\Xi }$   does not 
satisfy the left inverse property since there is  an element $\xi \in \Xi$ 
such that   $\xi^{-1}$ is not contained in $\Xi $.

If the loop  $L_{\Xi }$ satisfies the right inverse property 
then  using    
\[ (u,M) \ast (-(B-M^{-\delta }AM)^{-1} M^{-\delta }(B-A) u, M^{-1})
=(0,I) \] we have 
\[ (u, X)=[(u,X) \ast (u',M)] \ast (-(B-M^{-\delta }AM)^{-1} M^{-\delta }(B-A) u', M^{-1}) =  (u+ \] 
\[ [(B-(M X)^{ \delta }A (M X)^{-1})^{-1} X^{\delta }(B-M^{\delta } A 
M^{-1})- 
(B- X^{\delta }A X^{-1})^{-1} X^{\delta }(B-A)] u', X) \] 
for all $u, u' \in \mathbb K^n$ and $M,X  \in \Gamma_0$. 
For $X=M^{-1}$ this yields

\bigskip
\noindent
$(3)$ \centerline{$(B-A)^{-1}M^{-\delta }(B-M^{ \delta }A M^{-1})= (B-M^{-\delta }A M)^{-1} M^{-\delta }(B-A) $} 

\bigskip
\noindent
and for $X=M$ we obtain 

\bigskip
\noindent
$(4)$ \centerline{$ (B-M^{2 \delta }A M^{-2})^{-1}M^{\delta }(B-M^{\delta } A M^{-1})= 
(B- M^{\delta }A M^{-1})^{-1} M^{\delta }(B-A) $}

\bigskip
\noindent 
for all $M \in \Gamma_0$. 
Taking the inverses of both sides of (3) we get  

\bigskip
\noindent 
(5) \centerline{ $(B-M^{ \delta }A M^{-1})^{-1} M^{\delta } (B-A)=  (B-A)^{-1} M^{\delta } (B-M^{-\delta } A M)$. }

\bigskip
\noindent
Since the left side of (5) is equal to the right side of (4) we have 
for a proper loop $L_{\Xi }$ 
the contradiction $(B-A)=(B-M^{2 \delta } A M^{-2})$  for  all 
$M \in \Gamma_0$.

Let $T_B'$ be the minimal normal subgroup of $\Gamma $ containing $T_B$ and 
let $\Delta '$ be the semidirect product $T_B' \rtimes \Gamma _0'$. 
If 
$g(u_0,B u_0,I)$ is  an element of 
\[ T_B \cap \{ g( M_0 u, M_0^{\delta }B u, I); \ u \in \mathbb K^n \}=T_B \cap g(0,0,M_0) T_B g(0,0,M_0)^{-1} \]
then one has    
$u_0=M_0 u$ and $B u_0= M_0^{\delta} B u$ or $B M_0 u=M_0^{\delta } B u$. 
But in this case  the matrix  
$ M_0^{-1}B^{-1} M_0^{\delta } B$ 
would have an  eigenvalue $1$. This contradiction yields  that 
\[ T_B'=\{ g(u,v,I); \  u,v \in \mathbb K^n \}= T_B \times [g(0,0,M_0) T_B g(0,0,M_0)^{-1}] 
\]
and one has $\Delta '=\Gamma $.

The maximal normal subgroup $\Theta $ of $\Gamma $ contained in $H$ consists of the 
matrices $g(v, Av, I)$, where $v$ is an element of a subspace $V$  such that 
for all $v \in V$ and $M \in \Gamma _0$ one has 
\[ g(0,0,M) g(v, Av,I) g(0,0,M^{-1})=g(M v, M^{\delta } A v,I)=g(v', Av',I). \]
This is equivalent to $v'=M v$ and $ M^{\delta } A v= A M v$ for all 
$M \in \Gamma _0 $ and  $v \in V$. Hence  for the restrictions  of $M^{\delta } A$ and $A M$ to $V$ we have  $ M^{\delta } A|_V=  A M|_ V$ or 
equivalently $ M^{\delta } |_V =  A M A^{-1}|_ V$. 
According to Prop. 1.13 in \cite{loops}  the group generated by the left 
translations of the loop $L_{\Xi }$ is the group 
$\Delta=\Delta '/ (\Theta \cap \Delta ')$ and the  stabilizer of the identity 
$e \in L_{\Xi }$ is the group $ H / (\Theta  \cap \Delta ')$.

If $\mathbb K$ is a topological field, respectively  the field of real or 
complex  numbers then  $\Gamma $ is a 
topological group, respectively a  Lie group,   
and the section $\sigma $ is continuous, respectively analytic. Then  the multiplication 
of $L_{\Xi }$ as well as the left divison  
$(a,b) \mapsto a \backslash b: L_{\Xi } \times L_{\Xi } \to L_{\Xi }$ are 
continuous, respectively analytic. Looking at the solution of the equation 
$(1)$ we see that also the right divison  
$(a,b) \mapsto a / b: L_{\Xi } \times L_{\Xi } \to L_{\Xi }$ is 
continuous, respectively analytic.  
\qed

\bigskip
\noindent
The group $ \Gamma $ may be regarded as a subgroup  of the group of 
affinities of 
the $2n$-dimensional affine space ${\cal A}_{2n}$ acting on the set 
$\{ (1,x,y); \ x,y \in \mathbb K^n \}$    by 
\[ g(u,v,M) (1,x,y)=(1, u+M x, v+M^{\delta } y). \]
Then  $ \Gamma _0 '= \{ g(0,0,M); \ M \in  \Gamma _0 \}$ is  
the stabilizer of the point $(1,0,0) \in  {\cal A}_{2n}$ in  $\Gamma $ 
and  $\{g(u,v,I);\ u,v \in \mathbb K^n \}$ is  the 
translation group  of  ${\cal A}_{2n}$.
The $(2n-1)$-dimensional projective space  
\[ E=\{ \mathbb K^*(0,x,y);\ x,y \in  \mathbb K^n, (x,y) \neq (0,0), 
\mathbb K^*=\mathbb K \backslash \{ 0 \}  \} \] is the hyperplane at 
infinity of the affine space  ${\cal A}_{2n}$. The group $\Gamma $ acts on 
$E$ by $g(u,v,M)(0,x,y)=(0,u+M x, v+M^{\delta } y)$.

The group  $\Gamma _0 ' $ leaves the subspace 
${\cal Q}_1=\{ (1,x ,0);\ x \in  \mathbb K^n \}$ as well as the subspace 
${\cal Q}_2=\{ (1,0 ,x);\ x \in  \mathbb K^n \}$ invariant. 
We call a subspace ${\cal Q}_C$  of the form 
$\{ (1,x,Cx);\ x \in \mathbb K^n \}$ with $C \in GL(n, \mathbb K)$ an 
$n$-dimensional transversal subspace with  respect to $\Gamma _0' $. 
Any transversal subspace ${\cal Q}_C$ intersects ${\cal Q}_1$ and  
${\cal Q}_2$  only in the point $(1,0 ,0)$. The projective subspace 
 ${\cal Q}_C^*= \{ \mathbb K^* (0,x, C x);\ x \in \mathbb K^n \}$ of $E$ 
may be seen as the trace of ${\cal Q}_C$ in $E$.

\bigskip
\noindent
In the next Lemma we give necessary and sufficient  conditions  for 
the existence of 
regular orbits  for the group  
 $ \Gamma _0 '$ in the set  ${\cal T}$ of $n$-dimensional transversal 
subspaces of ${\cal A}_{2n}$. These conditions are needed for 
a geometric representation of loops $L_{\Xi }$ within the affine space 
${\cal A}_{2n}$.

\begin{Lemma} 
The group  $ \Gamma _0 '$ has in the set  ${\cal T}$ of $n$-dimensional affine 
transversal subspaces 
of  ${\cal A}_{2n}$  a regular orbit 
${\cal O}= \{ \varphi ({\cal Q}_A);\ \varphi \in \Gamma_0' \}$ if and only 
if one of the following equivalent conditions is satisfied:

\smallskip
\noindent
(i) There exists an inner automorphism $\alpha $ of $GL(n, \mathbb K)$  given by 
$X \mapsto A X A^{-1}$ with $A \in GL(n, \mathbb K)$ such that 
$M^{\delta } M^{-\alpha } \neq I$ for all 
$M \in \Gamma _0 \backslash \{ I \}$. 

\smallskip
\noindent
(ii) There exists an orbit  
${\widehat {\cal O}}= \{ \varphi  ({\cal Q}_A);\ \varphi \in \Gamma \}$  
such that the stabilizer of ${\cal Q}_A$ in $\Gamma $ is the group 
$H=\{ g(u,A u,I);\ u \in \mathbb K^n \}$.    
\end{Lemma}   
\Pro  
Let  ${\cal Q}_A=\{ (1,x,Ax);\ x \in \mathbb K^n \}$ be an $n$-dimensional transversal subspace of  ${\cal A}_{2n}$.  The orbit  ${\cal O}$ 
containing ${\cal Q}_A$ under $\Gamma_0'$ consists of all 
subspaces 
\[ \{ (1, M x, M^{\delta } A x); \ x \in  \mathbb K^n \}=\{ (1,x,M^{\delta } A M^{-1}x); \ x \in  \mathbb K^n \}, \] 
where $M$ varies over the elements of  $ \Gamma _0 $.   
The orbit  ${\cal O}$ is a regular orbit of  $ \Gamma _0 ' $ if and only if 
$A \neq M^{\delta } A M^{-1}$ or $I \neq  M^{\delta } M ^{- \alpha }$  for 
all 
$M \in  \Gamma _0 \backslash \{ I \}$, where $\alpha $ is the map  
$X \mapsto A X A^{-1}: GL(n, \mathbb K) \to GL(n,\mathbb K)$.

The stabilizer $\Gamma_{{\cal Q}_A}$ of ${\cal Q}_A$ in $\Gamma $ is the 
group 
$H=\{ g(u,Au, I);\ u \in \mathbb K^n \}$ if and 
only if 
the relation $g(u,v,M) (1,x,Ax)=(1,y,Ay)$ for all $x \in \mathbb K^n$ and suitable 
$y \in \mathbb K^n$ holds. Since $g(u,v,M) (1,x,Ax)=(1,u+M x, v+M^{\delta } A x)$ we 
obtain for $x=0$ that $v=A u$. Hence $H \le \Gamma_{{\cal Q}_A}$.  
Moreover, one has  $M^{\delta }A  x=AM x$ for all 
$x \in \mathbb K^n$. 
Therefore $H$ is the stabilizer of ${\cal Q}_A$ in $\Gamma $ if 
and only if 
for each $I \neq M \in \Gamma _0$ there is an $0 \neq x_0 \in \mathbb K^n $ 
such that $M^{\delta }A x_0 \neq  AM x_0$. This is equivalent to 
$M^{\delta } A M^{-1} \neq A$ which is the condition $(i)$.  
It follows that the conditions $(i)$ and $(ii)$  are equivalent for the existence of a 
regular orbit of $\Gamma_0'$ in the set ${\cal T}$. 
\qed

\bigskip
\noindent 
Using the geometric  interpretation for $\Gamma $ we prove that 
the loops  $L_{\Xi }$ have realizations  
on the orbit   ${\widehat {\cal O}}$ if the conditions 
of  Lemma 2 are satisfied  (cf. \cite{figula}, Theorem 1, p. 153).

\begin{Theo}
Let $L_{\Xi }$ be the loop determined by the matrices  
$A, B \in GL(n, \mathbb K)$, the group $\Gamma_0$ and the endomorphism 
$\delta :\Gamma _0 \to \Gamma _0$.  Let 
$H=\{g(u,Au,I);\ u \in \mathbb K^n \}$  be the stabilizer of 
${\cal Q}_A$ in $\Gamma $.   Then 
$L_{\Xi }$ is isomorphic to a  loop $L_{\Xi}'$ 
the elements of which are elements  of the orbit 
${\widehat {\cal O}}=\{ \psi ({\cal Q}_A); \ \psi \in \Xi \}=\{ \varphi ({\cal Q}_A);\ \varphi \in \Gamma \}$.   
The loop  $L_{\Xi}'$ has ${\cal Q}_A$ as the identity and the 
multiplication of  $L_{\Xi}'$ 
is defined by 
\[ X \circ Y=\tau _{{\cal Q}_A,X} (Y) \quad \hbox{for  all} \quad X,Y \in {\widehat {\cal O}}, \]
where $\tau _{{\cal Q}_A,X}$ is the unique element of $\Xi $  mapping 
${\cal Q}_A$ onto $X$.

The group $\Gamma _0' $ acts sharply transitively on the traces of elements 
of 
${\cal O}$ in the hyperplane $E$ at infinity, and  
the subspace $\{ (1,x, Bx); \ x \in \mathbb K^n \}$ 
  intersects any subspace in  ${\cal O}$ in 
precisely one  point.   
\end{Theo}
\Pro 
Since the subgroup $H$ is the stabilizer of ${\cal Q}_A$ in $\Gamma $ and 
the set $\Xi $ acts sharply transitively on the cosets $g(0,v,M) H$ the 
loop $L_{\Xi }$ 
is isomorphic to the loop $L_{\Xi }'$ defined on ${\widehat {\cal O}}$, with 
${\cal Q}_A$ 
as identity and with the multiplication defined  in the assertion. 

According to Lemma 2 the group $\Gamma_0'$ acts sharply transitively on the 
orbit ${\cal O}=\{ \varphi  ({\cal Q}_A);\ \varphi  \in \Gamma _0'\}$ and 
hence  also sharply transitively on the set of traces of elements of  
$\cal{O}$ in the hyperplane $E$.    

The subspace $ \{ (1,x,B x);\ x \in \mathbb K^n \}$  intersects any 
element  of   ${\cal O}$ 
only in the point  $(1,0,0)$  if and only if $(1,x,Bx) \neq (1,x, M^{\delta }A M^{-1}x) $ for all $M \in   \Gamma _0 \backslash \{ I \}$ and  $x \neq 0$. By  Theorem 1   
no of the matrices  $M^{-1} B^{-1} M^{\delta } A$ with  $M \in   \Gamma _0 $  has an  eigenvalue $1$. Hence   
$BM(I-M^{-1} B^{-1} M^{\delta } A )=$ $BM-M^{\delta } A =$
$(B-  M^{\delta } A M^{-1})M =B[(I- B^{-1}M^{\delta }A M^{-1})M]$, 
and the last claim of the Theorem  follows. 
\qed

\bigskip
\noindent
\centerline{\bf 3. Applications}

\bigskip
\noindent
In this Section we give concrete  examples for matrices $A,B \in GL(n, \mathbb K)$, groups $\Gamma _0$   and endomorphisms 
$\delta :\Gamma _0 \to \Gamma _0$ such that the loop 
 $L_{\Xi }=L_{A,B, \Gamma_0, \delta }$ exists.  To archieve this  goal we 
have in particular to show that no matrix  
 $M^{-1} B^{-1} M^{\delta } A$ has an  eigenvalue $1$ for all 
$M \in \Gamma_0$. 

\bigskip
\noindent
{\bf 3.1}  Let $\mathbb K$ be a commutative field and  let  $\Gamma _0$ be a 
subgroup of $GL(n, \mathbb K)$.  
   
\smallskip
\noindent
a) We assume that  the group $\Gamma_0$ is not commutative and that 
 $\delta $ is the inner automorphism $X \mapsto C^{-1} X C$ of $\Gamma _0$ different from the identity. We choose $A=C^{-1}$ and $B \in GL(n, \mathbb K)$ 
 such  that $B^{-1} C^{-1}$ does not centralize $\Gamma_0$ and  has no 
eigenvalue $1$. Then the loop 
$L_{A,B, \Gamma_0, \delta  }$ is a proper loop. But  because of $M^{\delta } C^{-1} M^{-1} C=I$ for all $M \in \Gamma_0$  this  loop  has  no  geometric 
realization in sense of Theorem 3.

\smallskip
\noindent
b) We suppose that $\delta $ is the identity, the matrix $A$ centralizes 
 $\Gamma_0$ but  the matrix $B$ does  not centralize  $\Gamma_0$. 
If  $B^{-1} A$ has no  eigenvalue  $1$ then  the eigenvalues of 
the  matrices  $M^{-1} B^{-1} A M$ are also 
different from $1$ for all  
$M \in \Gamma_0$. Hence the proper loop 
$L_{A,B, \Gamma_0, id }$ exists.

\smallskip
\noindent 
c) Let $P_m$ be an $(m \times m)$-matrix and let $Q_{n-m}$ be an 
$(n-m \times n-m)$-matrix. 
We denote by $P_m \oplus Q_{n-m}$ the matrix $\left (\begin{array}{cc}
P_m & 0 \\
0 & Q_{n-m} \end{array} \right)$. 
Let  $A=diag (a, \cdots ,a) \oplus A'$ with 
$diag (a, \cdots ,a) \in GL(m,\mathbb K)$ and $A' \in GL(n-m, \mathbb K)$, 
let  $B=B' \oplus diag (b, \cdots ,b)$ with 
$diag (b, \cdots ,b) \in GL(n-m,\mathbb K)$ and $B' \in GL(m,\mathbb K)$. 
Let $\widehat{\Gamma_0 }  \neq I$  be a subgroup of 
$GL(m,\mathbb K)$ and 
$\Gamma_0=\{ M \oplus I_{n-m}; M \in  \widehat{\Gamma_0 } \}$, where 
$I_{n-m}$ is the 
identity  of $GL(n-m, \mathbb K)$. 
If  $\delta $ is the identity  the loop $L_{A,B,\Gamma_0, id}$ 
exists if 
and only if  
$A'$ has no eigenvalue $b$,  
the matrix  $B'$ has no eigenvalue $a$ and $a \neq b$ holds.  Futhermore,  
$L_{A,B,\Gamma_0, id}$ is a proper loop if $B'^{-1}$ does not centralize  
$\widehat{\Gamma_0 }$. 

\smallskip
\noindent
The loops  $L_{A,B,\Gamma_0, id}$ treated in b) and c) 
 have  no  geometric realizations in sense of Theorem 3 since  
 $A$  centralizes $\Gamma_0$.

\bigskip
\noindent
{\bf 3.2}  Let $\mathbb K$ be a commutative field. Let $\Gamma _0$ be a 
non-abelian  
subgroup  of the group of upper triagonal matrices 
whose entries are elements of  $\mathbb K$.  Let 
$\chi _i$, $i=1, \cdots,n,$ be the map which 
assigns to the matrix $M=(m_{ij})$ the element $m_{ii}$. Then   
$\chi _i $ is a   homomorphism  from $\Gamma _0$ into the 
multiplicative group $\mathbb K^*$  of $\mathbb K$.     
Let $\delta =0$ be the  
endomorphism which maps  $\Gamma _0$ onto  $I$. 
If $B^{-1}A=\hbox{diag}\ (t_1,t_2, \cdots,t_n)$  
such that $\Phi _i=\{ \chi_i(M);\ M \in \Gamma_0 \}$ is a proper 
subgroup of  $\mathbb K^*$ 
and  $t_i \not \in \Phi _i$ for 
$i=1, \cdots ,n$, then  $L_{A,B, \Gamma_0, 0 }$ is a proper loop  
and has a 
geometric realization on the set 
$\{ \varphi  ({\cal Q}_A);\ \varphi \in \Gamma \}$ 
 (Theorem 3 and Lemma 2).

\bigskip
\noindent
{\bf 3.3} Let $\Gamma_0$ be the 
group $SU_2(\mathbb C)=\left \{ \left (\begin{array}{rr} 
z & w \\
-\bar{w} & \bar{z} \end{array} \right); z,w \in \mathbb C,z \bar{z}+ w \bar{w}=1 \right \}$. 
Every element of  $\Gamma_0$ has 
eigenvalues $\{ e^{i \Theta }, e^{-i \Theta} \}$ with 
$0 \le \Theta < 2 \pi$. Let $\delta $ be the inner automorphism  
$X \mapsto C^{-1} X C$ of $\Gamma _0$, let $A$ and $B$ be elements of 
$\Gamma_0$. We assume that $B^{-1} C^{-1}$ has an eigenvalue 
$e^{i \Theta _1}$ 
 and $CA$ has an eigenvalue  $e^{i \Theta _2}$ with 
$0< \Theta _i $, $i=1,2$,  $ \Theta _1 > \Theta _2$ and 
$ \Theta _1 +\Theta _2 \le \pi$. Since  for any 
$M \in \Gamma_0$ the matrix $M^{-1} B^{-1} C^{-1} M$ has an  eigenvalue  
$e^{i \Theta _1}$  we 
have  (see  \cite{weitsman}, Prop. 3.1, p. 601)  the inequalities   
\newline
\centerline{$\cos ( \Theta _1 +\Theta _2 ) \le \cos  \Theta _3 \le \cos ( \Theta _1 -\Theta _2 )$, }
where $e^{i \Theta _3}$ is  an eigenvalue of the matrix 
$[M^{-1} B^{-1} C^{-1} M] C A$.  It follows that no matrix 
$[M^{-1} B^{-1} C^{-1} M] C A$ has the eigenvalue $1$ and  the 
differentiable  proper loop 
 $L_{A,B, \Gamma_0, \delta }$ exists. Since $-I$ is contained in the centre 
of $SU_2(\mathbb C)$ the loop $L_{A,B, \Gamma_0 , \delta }$ has no geometric 
realization.

\bigskip
\noindent
{\bf 3.4} (i) Let $\Gamma_0 \neq I$ be a  compact  subgroup of 
$GL(n, \mathbb R)$.

\smallskip
\noindent
(ii) Let $\mathbb K$ be a commutative field with an exponential 
(ultrametric) valuation $v: \mathbb K \to \mathbb R \cup \{ \infty \}$ 
(cf. \cite{endler}, p. 20, \cite{serre}, p. 65) and let ${\mathcal A}$ 
be the corresponding valuation ring.  The field  $\mathbb K$ is 
 a topological field with respect to $v$ (\cite{endler}, p. 2) and 
 $GL(n, \mathbb K)$  
carries the topology induced by the topology of  $\mathbb K$. 
Let $M$ be a matrix which topologically generates a compact subgroup 
$\Upsilon $ of $GL(n, \mathbb K)$. The matrix $M$  is  conjugate to an upper triangular matrix in 
$GL(n, \mathbb L)$, where $\mathbb L $ 
is a finite algebraic extension of  $\mathbb K$. Let ${\hat v}$ be an 
extension of the valuation of $\mathbb K$ to 
$\mathbb L$ and let ${\hat {\mathcal A}}$ be the corresponding valuation 
ring.  Assume that $\lambda \neq 0$ is an eigenvalue of  $M$. Since 
$\Upsilon $ is compact there exist natural numbers 
$n_i$ 
with $\lim \limits_{i \to \infty} n_i=\infty $ such that $M^{n_i}$ 
converges to $S$. The matrix $S$ has 
$\lim \limits_{i \to \infty} \lambda ^{n_i}$ as an eigenvalue. Because of 
${\hat v}(\lim \limits_{i \to \infty} \lambda ^{n_i})=\lim \limits_{i \to \infty} 
{\hat v}( \lambda ^{n_i})= {\hat v}(\lambda ) 
 \lim \limits_{i \to \infty} n_i$ and  ${\hat v}(x)=\infty $ if and only 
if $x=0$ 
it follows that ${\hat v}(\lambda )=0$. This means that $\lambda $ is a 
unit in the valuation ring ${\hat {\mathcal A}}$.  

Let 
$\Gamma _0 \neq I$ be a closed subgroup of the group $GL(n, {\mathcal A})$. 
According to \cite{serre}, p. 104,  the group $GL(n, {\mathcal A})$ is 
compact and  hence $\Gamma _0$    
 is also  compact.  

\smallskip
\noindent
Let $A=diag (a, \cdots ,a)$ and 
$B=diag (b, \cdots ,b)$ 
be diagonal matrices in the centre of $GL(n,\mathbb F)$. If  $\mathbb F$ 
is the field of real numbers then we suppose that $| b^{-1} a| \neq 1$. If 
$\mathbb F$ has an  exponential valuation then we assume that 
$v( a b^{-1}) \neq 0$.  Then any matrix 
$M^{-1} B^{-1} M^{\delta } A=M^{-1} M^{\delta } B^{-1} A$ with 
$M \in \Gamma _0$ has no eigenvalue $1$ (cf. \cite{zurmuhl},  p. 288  and 
Satz 8, p. 196).

If $\delta $ is not the identity then  $L_{A,B, \Gamma_0, \delta }$ is a 
proper loop. It has  a geometric realization if and only if 
$M^{\delta } \neq M$ 
for all $M \in \Gamma_0 \backslash \{ I \}$. This is for instance the case 
if there exists a 
natural number $k$  such that $\delta ^k =0$. (If for a matrix 
$I \neq M \in \Gamma_0$ one has $M^{{\delta } ^k}=I$  then 
$M^{\delta } \neq M$.) Let 
$\Gamma_0 =H_1 \times H_2 \times \cdots \times H_k \le GL(n, \mathbb F)$  
such that $H_i$ is isomorphic to the same  compact Lie group, respectively  
closed  subgroup of $GL(n, {\mathcal A})$ for 
$i \in \{ 1, \cdots , k \}$.  Then  
$\delta: (x_1, \cdots , x_k) \mapsto (1,x_1, \cdots ,x_{k-1})$ is an 
endomorphism of $\Gamma_0$ and  $\delta ^k=0$.

\bigskip
\noindent
{\bf 3.5} 
(i) Let $\mathbb K$ be a commutative  field and let $\mathbb F$ be a 
subfield of  $\mathbb K$. Moreover,   let  
\[ \Gamma_0=\left \{ M=\left ( \begin{array}{cc}
m_{11} & m_{12} \\
m_{21} & m_{22} \end{array} \right);\ det (M)=1 \right \} \] 
  be the group $SL(2,\mathbb F)$. 
We assume that $\delta $ is the identity. 
\newline
Let $B=\left ( \begin{array}{rr}
1 & -t \\
0 & 1 \end{array} \right)$ and $A=\left ( \begin{array}{cc}
a & 0 \\
c & a \end{array} \right)$ be matrices with $t,c \in \mathbb F$ such that 
$ct$ is  
a square  in $\mathbb F$ and $a \in \mathbb K \backslash \mathbb F$.  
The trace $t(S_M)$ of the matrix  $S_M=M^{-1} B^{-1} M A$, where 
$M \in \Gamma_0$, has the value 
$ctm_{22}^2+2 a= \lambda _1 + \lambda_2$, where  $\lambda_i$, 
$i=1,2$, are the eigenvalues of $S_M$. If  $\lambda_1=1$ then  
$det (S_M)=a^2= \lambda_2 $ and the equation 
$a^2- 2 a +1 -ct m_{22}^2=0$ yields that  $a=1 \pm m_{22} \sqrt{ct} \in 
\mathbb F$, which is a contradiction. 
Hence the matrix  $S_M$ has no 
eigenvalue $1$ for all $M \in SL(2, \mathbb F)$. 

\smallskip
\noindent
(ii) Let $\mathbb K$ be a formally real field (i.e. the sums of squares are  
squares in $\mathbb K$, but $-1$ is not a square in $\mathbb K$) and let 
$\Gamma _0$ be 
the group $SL(2,\mathbb K)$. Moreover,  let $B=\left ( \begin{array}{rr}
a & -b \\
b & a \end{array} \right)$ and $A=\left ( \begin{array}{cc}
1 & t \\
0 & 1 \end{array} \right)$ be matrices in  $SL(2, \mathbb K)$ 
 such that  $1-a$ is a square but  $-tb$ is not a square in $\mathbb K$. 
The trace  $t(S_M)$ of the matrix  $S_M=M^{-1} B^{-1} M A$ has the value 
$-tb(m_{11}^2+m_{12}^2)+ 2a$. As $2(1-a)$ is a square, but  
$-tb(m_{11}^2+m_{12}^2)$ is 
not a square we have $-tb(m_{11}^2+m_{12}^2) \neq  2(1-a)$ for all  
$M \in SL(2,\mathbb K)$. 

Any loop $L_{A,B, \Gamma_0, id }$ of (i) as well as of (ii) is a 
proper loop
because the condition $(II)$ in Theorem $1$ is not satisfied.  
But it  has  no 
geometric realization since  $-I \in \Gamma _0$.

\bigskip
\noindent
{\bf 3.6} Let $\mathbb K$ be a commutative field and let $\Gamma_0$ be the 
group  
\[ \left \{ M:=g(m,n)=\left ( \begin{array}{ll}
m & n \\
0 & m^{-1} \end{array} \right);\  n \in \mathbb K, m \in \Omega, \right \} \]
where $\Omega $ is a subgroup  of $\mathbb K^*$ which does not 
contain $-1 \neq 1$. (If the characteristic of $\mathbb K$ is $2$ then we 
may take $\Omega =\mathbb K ^*$.) Let 
$\delta $ be the mapping $g(m,n) \mapsto g(m,d n)$ with 
$d \in \mathbb K$. If $d \neq 0$ then $\delta $ is an 
automorphism of  $\Gamma_0$, if $d=0$ then $\delta $ is a proper 
endomorphism 
of $\Gamma_0$. Let $B=\left ( \begin{array}{cc}
a^{-1} & - b \\
0 & a \end{array} \right)$ and $A=\left ( \begin{array}{cc}
r & s \\
0 & r^{-1} \end{array} \right)$ be  matrices in $SL(2, \mathbb K)$. 
We 
assume that  $a r \neq 1$. 
The trace  $t(S_M)$ of the matrix  $S_M=M^{-1} B^{-1} M^{\delta } A$ has 
the value $ar+ a^{-1} r^{-1} $. 
Hence  $t(S_M) \neq  2$ for all $M \in \Gamma_0$. 

The loop  $L_{A,B, \Gamma_0, \delta }$ 
is a proper loop if and only if  $a^2 d \neq 1$ or $b \neq 0$ because of the condition $(II)$ in Theorem $1$. 
 Moreover, it has  a geometric realization on the set 
$\{ \varphi  ({\cal Q}_A);\ \varphi \in \Gamma \}$  if and only if  
$M^{\delta } A \neq  A M$ for all 
$M \in \Gamma_0 \backslash \{ I \}$ or  
  $0 \neq sr (1-m^2)+nm (r^2-d)$ 
for all $(m,n) \neq (1,0)$ and $m \neq 0$.  This is for instance the 
case  if  $r^2=d$ and $s \neq 0$.

\bigskip
\noindent
{\bf 3.7} Let $\Gamma_0$ be the group of matrices 
$M:=g(m)=\left ( \begin{array}{ll}
m & 0 \\
0 & m^{-1} \end{array} \right)$, $0 < m \in \mathbb R$, and let $\delta $ be 
the automorphism $g(m) \mapsto g(m^c)$ 
with $0 \neq c \in \mathbb R$. Let 
$B^{-1}=\left ( \begin{array}{cc}
k & l \\
n & s \end{array} \right)$ and $A=\left ( \begin{array}{cc}
p & q \\
r & v \end{array} \right)$ be matrices of $SL(2, \mathbb R)$. For  the trace 
$t(S_M)$ of the matrix 
$S_M=M^{-1} B^{-1} M^{\delta } A$ one has  
$d(m^{c-1}+m^{-(c-1)})+lr m^{-(c+1)}+nq m^{(c+1)} $, where $d:=kp= sv$. 
If $d>1$,  $lr \ge 0$ and  $nq \ge 0$ or $d<-1$, 
$lr \le 0$ and  $nq \le 0$ then $t(S_M) \neq 2$ for all $M \in \Gamma_0$ 
and  $S_M$ has  no eigenvalue $1$. 

For $c \neq 1$ the loop  $L_{A,B, \Gamma_0, \delta }$ is always a proper 
loop which has a  geometric realization on the 
set $\{ \varphi  ({\cal Q}_A);\ \varphi \in \Gamma \}$. If $c=1$ then the 
loop  $L_{A,B, \Gamma_0, \delta }$ is a proper loop 
 if and only if  $l \neq 0$ or $n \neq 0$ and  has  a geometric realization 
 precisely if $q \neq 0$ or $r \neq 0$ (cf. condition 
$(II)$ in Theorem $1$ and Lemma $2$).

\bigskip
\noindent
{\bf 3.8} Let $\Gamma_0$ be the group consisting of the matrices 
\[M:=g(\varphi )=\left ( \begin{array}{rr}
\cos{\varphi } & \sin{\varphi } \\
-\sin{\varphi } & \cos{\varphi } \end{array} \right),\  \varphi \in [0, 2 \pi),  \]
and let $\delta $ be the mapping  $g(\varphi ) \mapsto g(n \varphi )$
with $n \in \mathbb Z$. If  $n \not \in \{-1, 1\}$ then  
$\delta $ is 
a proper  endomorphism of $\Gamma_0$, if $n \in \{-1,1 \}$  then 
$\delta $ is 
an automorphism of  $\Gamma_0$.  Let $B^{-1}=\left ( \begin{array}{rr}
a & b \\
c & -a \end{array} \right)$, $a,b,c \in \mathbb R$, $-a^2-bc=1$, 
and let  $A=\left ( \begin{array}{rr}
\frac{1}{\sqrt{2}} & \frac{1}{\sqrt{2}}  \\
- \frac{1}{\sqrt{2}} &  \frac{1}{\sqrt{2}} \end{array} \right)$. 
If $b>0$, $c<0$ and $-\sqrt{2} < c-b <0$ then for the trace $t(S_M)$ of 
$S_M=M^{-1} B^{-1} M^{\delta } A$ one has  
\newline
\centerline{ $t(S_M)=\displaystyle \frac{c-b}{\sqrt{2}}[ \cos{((n-1)\varphi )}+\sin{((n-1)\varphi )}] < \frac{2(b-c)}{\sqrt{2}}<2$. }
 Hence the eigenvalues of $S_M$ 
for all $M  \in \Gamma_0$ are different from $1$. 

Any 
loop $L_{A,B, \Gamma_0, \delta }$ is a proper loop since the condition 
$(II)$ in Theorem $1$ is not satisfied.   
Moreover,  it   
 has a geometric realization on the set 
$\{ \varphi  ({\cal Q}_A);\ \varphi \in \Gamma \}$  precisely if  
$M^{\delta } A \neq  A M$  for all  $M \in \Gamma_0 
\backslash \{ I \}$. This is the case if and only if $n \neq 1$.

\bigskip
\noindent
{\bf 3.9} Let $a,b$ be  elements of the multiplicative 
group  $\mathbb K^*$ of a commutative field  $\mathbb K$.  A  proper loop 
$L_{a,b,\Gamma _0, \delta}$ exists if one of  the following conditions holds:

\smallskip
a) $\Gamma_0 \neq 1$ is a proper subgroup of  $\mathbb K^*$, the endomorphism 
$\delta $ is fixed point-free (i.e. $x^{\delta } \neq x$ for all 
$x \in \Gamma _0 \backslash  \{ 1 \}$) and $b^{-1}a \notin \Gamma _0$. This 
is for instance the case if $\Gamma_0$ does not contain $-1 \neq 1$ and 
$\delta $ is the automorphism $x \mapsto x^{-1}$. If $\mathbb K$ has the 
characteristic $2$ then any proper subgroup of  $\mathbb K^*$ is suitable 
 as $\Gamma_0$.

\smallskip
b) $\mathbb K = \mathbb R$, $\Gamma_0$ is the multiplicative group 
$\mathbb R^*$, $b^{-1}a <0$  and $\delta:x \mapsto x^{2k+1}$ 
for an integer $k \neq 0$. 

\smallskip
c) $\mathbb K=\mathbb C$ and $\Gamma_0$ is the multiplicative group  
$\mathbb C^*$. If $b^{-1}a$ is a real number then let 
$\delta $ be the endomorphism $r e^{it} \mapsto r e^{-it}$. If $b^{-1}a$ is not real  then let 
$\delta $ be the endomorphism $r e^{it} \mapsto r^{-1} e^{it}$.  

The loop $L_{a,b,\Gamma_0, \delta }$ has a geometric realization on the set 
$\{ \varphi  ({\cal Q}_A);\ \varphi \in \Gamma \}$  in the case (a), but 
 no  geometric realization in the cases (b) or (c).

\bigskip
\noindent
\centerline{\bf 4.  Groups  generated by  right translations }

\bigskip
\noindent
Let  $L_{\Xi}$ be a loop realized on the semidirect product 
 $\mathbb K^n \rtimes \Gamma _0$ and let 
 $\varphi $ be the epimorphism $L_{\Xi} \to \Gamma_0$ mapping $(u,M)$ 
onto $M$. If  $\varrho _a$ is 
 a  right translation of $L_{\Xi}$ then    
$\varrho _{\varphi (a)}$ denotes the corresponding right translation of 
$\Gamma_0$ by $\varphi (a)$.  
Because of 
$\varrho _{\varphi (a b)}=\varrho _{\varphi (a) \varphi (b)}=
\varrho _{\varphi (a)} 
\varrho _{\varphi (b)}$, there is an epimorphism $\omega  $ from the group 
$\Sigma $ generated by the right translations of  $L_{\Xi}$ onto 
$\Gamma_0$, such that 
$\varrho _{(u,M)}=\varrho_M$. The kernel $N$ of $\omega $ consists of all 
elements of
 $\Sigma $ which leave any subset $P_M=\{ (u,M);\ u \in \mathbb K^n \}$ 
invariant. Since $\Sigma $ acts sharply transitively on the set 
$\{P_M;\ M \in \Gamma_0 \}$ the group $\Sigma $ is a semidirect product 
$N \rtimes \Gamma_0$ of $N$ by $\Gamma_0$.

\bigskip
\noindent
If  $\mathbb K=\mathbb R$ then  $\Sigma $ is a smooth group and the 
manifold $P_M$ is diffeomorphic 
to $\mathbb R^n$.  

\begin{Prop} The group $\Sigma $ topologically generated by all  right 
translations of the proper smooth  
loop  $L_{\Xi}$ is a  smooth group  which contains an 
infinite dimensional abelian subgroup ${\mathcal D}$. The subgroup 
${\mathcal D}$ leaves any manifold $P_M$  invariant.   
\end{Prop}

\Pro The right translation of  $L_{\Xi}$ by $(u_2,M_2)$ is the smooth map
\[ \varrho _{(u_2,M_2)}: (u_1,M_1) \mapsto (u_1+u_2^{\psi }, M_1 M_2), \]
where  $\psi :\mathbb R^n \to \mathbb R^n$ is the linear  map 
defined by (I) in Theorem 1. 
Hence  $\Sigma $ contains the subgroup $S=\{ \varrho _{(u,I)};\ u 
\in \mathbb R^n \}$ of the mappings 
\[ \varrho _{(u,I)}: (u_1,M_1) \mapsto (u_1+ [B-M_1^{\delta }A M_1^{-1}]^{-1}M_1^{\delta }
(B-A)u,M_1). \] 
The group $S$, which is contained in the normal subgroup $N$ of $\Sigma $,  
is  diffeomorphic to $\mathbb R^n$. The conjugate subgroup 
\[\Sigma _{M}=\varrho ^{-1}_{(0,M)} S \varrho _{(0,M)}=
\varrho _{(0,M^{-1})} S \varrho _{(0,M)} \]
 consists  of the mappings 
\[ (u_1,M_1) \mapsto (u_1+[B-(M_1M)^{\delta }A (M_1M)^{-1}]^{-1}(M_1 M)^{\delta }(B-A)u,M_1). \] 
The set of subgroups $\Sigma _{M}, M \in \Gamma_0$,  generates in the group 
$\Sigma $ an abelian subgroup $\mathcal D$, which is a real vector space 
contained in $N$.    
We assume that  $\mathcal D$ has  finite dimension. 
Let $0 \neq u \in \mathbb R^n$ be a  fixed vector.  Then there exist 
elements   
\[  ([B-(M_1M^{(i)})^{\delta }A (M_1M^{(i)})^{-1}]^{-1}(M_1 M^{(i)})^{\delta }(B-A)u, M_1) ,\ i=1, \cdots , m,  \] 
such that from matrix equation 

\bigskip
\noindent
$(6)$ \centerline{$\sum \limits_{i=1}^m\  \nu _i [B-(M^{(i)})^{\delta }A 
(M^{(i)})^{-1}]^{-1} ( M^{(i)})^{\delta }=0, \ \nu _i \in \mathbb R $,}

\bigskip
\noindent
it follows $\nu _i=0$ for all $i=1, \cdots ,m,$.  Moreover, 
for  any $M^* \in \Gamma_0$ there are   
real  numbers  
$\lambda _i,\ i=1, \cdots , m,$ satisfying the identity  

\bigskip
\noindent
$(7)$ \centerline{$\sum \limits_{i=1}^m\  \lambda _i [B-(M_1M^{(i)})^{\delta }A (M_1M^{(i)})^{-1}]^{-1} (M_1 M^{(i)})^{\delta }=$}

\bigskip
\noindent
\centerline{$[B-(M_1M^*)^{\delta }A (M_1M^*)^{-1}]^{-1}(M_1 M^*)^{\delta }$}  

\bigskip
\noindent
for all  $M_1 \in \Gamma_0$.  For $M^* \in \{M_1^{-1}, I,M_1 \}$ the 
equation $(7)$ yields 
\[ (B-A)^{-1}=\sum \limits_{i=1}^m\  \lambda _i [B-(M_1M^{(i)})^{\delta }A 
(M_1M^{(i)})^{-1}]^{-1} (M_1 M^{(i)})^{\delta } \]
\[ [B-M_1^{\delta }A M_1^{-1}]^{-1}M_1^{\delta }= \sum \limits_{i=1}^m\  \mu  _i [B-(M_1M^{(i)})^{\delta }A 
(M_1M^{(i)})^{-1}]^{-1} (M_1 M^{(i)})^{\delta } \]
and 
\[ [B-M_1^{2 \delta }A M_1^{-2}]^{-1}M_1^{2 \delta }= \sum \limits_{i=1}^m\  \nu  _i [B-(M_1M^{(i)})^{\delta }A 
(M_1M^{(i)})^{-1}]^{-1} (M_1 M^{(i)})^{\delta } \] for suitable $\lambda _i, 
\mu _i, \nu _i \in \mathbb R$. 
Putting in these equations $M_1=I$ and using  $(6)$ we obtain 
$\lambda _i=\mu _i= \nu _i$ for all $i=1, \cdots ,m,$ and 
\[  (B-A)^{-1}=[B-M_1^{\delta }A M_1^{-1}]^{-1}M_1^{\delta }= [B-M_1^{2 \delta }A M_1^{-2}]^{-1}M_1^{2 \delta }. \] 
For $\delta =0$ one has  $(B-A)^{-1}=[B-A M_1^{-1}]^{-1}$ or  
$(B-A)= [B-A M_1^{-1}]$
which gives  the contradiction $0=A (I-M_1^{-1})$ for all 
 $M_1 \in \Gamma_0$. 

If  $\delta \neq 0$ then we obtain $[B-M_1^{\delta }A M_1^{-1}]^{-1}M_1^{\delta }=$ $[B-M_1^{2 \delta }A M_1^{-2}]^{-1}M_1^{2 \delta }$
or $(I-M_1^{-2 \delta })B =A (I-M_1^{-1})$.  
Hence $M_1^{\delta }= M_1^{2 \delta }$ for all 
$M_1 \in \Gamma_0$ which  is  a contradiction.   
\qed

\bigskip
\noindent
\centerline{\bf 5. The Akivis algebra of the smooth loop $L_{A,B,\Gamma _0, \delta }$ }

\bigskip
\noindent
Let $\Gamma_0$ be a  Lie subgroup of positive dimension in  
$GL(n,\mathbb R)$. Then the Akivis algebra 
${\mathfrak a}_{L_{\Xi }}=({\mathfrak a}_{L_{\Xi }}, [.,.], \langle .,.,. 
\rangle ) $ of a smooth loop 
$L_{\Xi }=L_{A,B,\Gamma _0, \delta }$ can be obtained in the following way. 
Let $\exp  m$ be the exponential image of the element $m$ in the Lie algebra 
${\mathfrak m}$ of  $\Gamma_0$. Let 
\[ C_{i,j}= (x_i, \exp{m_i}) \ast (x_j, \exp{m_j})= \] 
\[ \big (x_i+(B-(\exp{m_i} \exp{m_j})^{\delta } A (\exp{m_i} \exp{m_j})^{-1})^{-1} \] 
\[ (\exp{m_i})^{\delta } (B-(\exp{m_j})^{\delta } A (\exp{m_j})^{-1}) x_j,
 \exp{m_i} \exp{m_j} \big ),  \] 
where $i,j \in \{ 1,2 \}$. 
One has 
\[ C_{1,2}/C_{2,1}= \big (I-[B-(\exp{m_1} \exp{m_2})^{\delta } A (\exp{m_1} \exp{m_2})^{-1}]^{-1} \] 
\[ (\exp{m_1} \exp{m_2} \exp{m_1}^{-1})^{\delta }
 (B- \exp{m_1}^{\delta } A \exp{m_1}^{-1}) \big ) x_1+ \] 
\[ [B-(\exp{m_1} \exp{m_2})^{\delta } A (\exp{m_1} \exp{m_2})^{-1}]^{-1} \]
\[ \{ \exp{m_1}^{\delta } (B- \exp{m_2}^{\delta } A  \exp{m_2}^{-1}) -  
 (\exp{m_1} \exp{m_2} \exp{m_1}^{-1} \exp{m_2}^{-1})^{\delta } \] 
\[ (B-(\exp{m_2} \exp{m_1})^{\delta } A (\exp{m_2} \exp{m_1})^{-1}) \} x_2,  \exp{m_1} \exp{m_2} \exp{m_1}^{-1} \exp{m_2}^{-1} \big). \]
Let 
\[ D_1=\big ( (x_1,\exp{m_1}) \ast (x_2,\exp{m_2}) \big ) \ast (x_3,\exp{m_3})= (x_1+ \]
\[[B-(\exp{m_1} \exp{m_2})^{\delta } A (\exp{m_1} 
\exp{m_2})^{-1}]^{-1} \exp{m_1}^{\delta } (B- \exp{m_2}^{\delta } A  \exp{m_2}^{-1})x_2\ + \]
\[ [B-(\exp{m_1} \exp{m_2} \exp{m_3})^{\delta } A (\exp{m_1} \exp{m_2} \exp{m_3})^{-1}]^{-1} \]
\[  (\exp{m_1} \exp{m_2})^{\delta } (B-  \exp{m_3}^{\delta } A (\exp{m_3})^{-1}) x_3,\  \exp{m_1} \exp{m_2} \exp{m_3}) \] 
and \[ D_2=(x_1,\exp{m_1}) \ast \big ( (x_2,\exp{m_2})) \ast (x_3,\exp{m_3}) \big ) = \]   \[ (x_1+[B-(\exp{m_1} \exp{m_2} \exp{m_3})^{\delta } A (\exp{m_1} \exp{m_2} \exp{m_3})^{-1}]^{-1} \]
\[ \exp{m_1}^{\delta } (B- (\exp{m_2} \exp{m_3})^{\delta } A (\exp{m_2} \exp{m_3})^{-1})x_2\ + \]
\[  [B-(\exp{m_1} \exp{m_2} \exp{m_3})^{\delta } A (\exp{m_1} \exp{m_2}  \exp{m_3})^{-1}]^{-1} (\exp{m_1} \exp{m_2})^{\delta } \] 
\[ (B-  \exp{m_3}^{\delta } A (\exp{m_3})^{-1}) x_3,\  \exp{m_1} \exp{m_2} \exp{m_3}). \]
Then one has 
\newline
\centerline{$D_1/D_2=$}
\[\big ( [B-(\exp{m_1} \exp{m_2})^{\delta } A (\exp{m_1} \exp{m_2})^{-1}]^{-1} \exp{m_1}^{\delta }  (B- \exp{m_2}^{\delta } A  \exp{m_2}^{-1})x_2 - \]
\[ [B-(\exp{m_1} \exp{m_2} \exp{m_3})^{\delta } A (\exp{m_1} \exp{m_2}  \exp{m_3})^{-1}]^{-1} \exp{m_1}^{\delta } \]
\[(B- (\exp{m_2} \exp{m_3})^{\delta } A (\exp{m_2} \exp{m_3})^{-1}) x_2,\  I). \]
To obtain the binary,  respectively the ternary operation of the Akivis algebra 
${\mathfrak a}_{L_{\Xi }}$, which is realized on the vector space 
$\mathbb R^n \oplus \mathfrak{m}$, we replace in 
$C_{1,2}/C_{2,1},$ respectively in  $D_1/D_2$
the elements $\exp {m_k}, k=1,2,$  by one parameter 
subgroups $\exp{t m_k}$,  the elements $x_k$  by one parameter 
subgroups  $t x_k$  and form the following limits:  
\[ \lim \limits_{t \to 0} \frac{1}{t^2} (C_{1,2}(t)/C_{2,1}(t))=:[(x_1,m_1),(x_2,m_2)], \]
\[ \lim \limits_{t \to 0} \frac{1}{t^3} (D_1(t)/D_2(t))=:  \langle (x_1,m_1), (x_2,m_2), (x_3,m_3) \rangle \]
(cf. \cite{hofmann}, Prop. 3.3, p. 323). Using often   the fact  
\[ \frac{d}{dt} (F(t))^{-1}=- (F(t))^{-1}  \frac{d}{dt} (F(t)) (F(t))^{-1} \]
we obtain by straightforward calculation that 

\bigskip
\noindent
(8) \centerline{$[(x_1,m_1),(x_2,m_2)]=$ }
\[ \big ( (B-A)^{-1} \{ (m_1^{\tilde {\delta }} B-A m_1)x_2+ (Am_2 - m_2^{\tilde {\delta }}B)x_1 \}, [m_1,m_2] \big ), \] 
as well as 
\[ \langle (x_1,m_1), (x_2,m_2), (x_3,m_3) \rangle = \]
\[ \big((B-A)^{-1} \big \{  (m_3^{\tilde {\delta }}A-Am_3)(B-A)^{-1} (Am_1-m_1^{\tilde {\delta }}B)-(m_3^{\tilde {\delta }}A m_1-Am_3 m_1) \big \} x_2,0 \big ), \] 
where in both cases $\tilde {\delta }$  is the endomorphism of 
$\mathfrak{m}$ corresponding to $\delta $.

A straightforward but  tedious calculation  shows that for the Akivis 
algebra 
${\mathfrak a}_{L_{\Xi }}$ the left as well as the right 
side of the Akivis identity equals to 

\bigskip
\noindent
$(9)$ \centerline{$\big( (B-A)^{-1} \big \{ (m_2^{\tilde {\delta }}A-Am_2)(B-A)^{-1} \big [Am_3-m_3^{\tilde {\delta }}B \big ]+ $}
\[ (Am_3-m_3^{\tilde {\delta }}A )(B-A)^{-1}[Am_2-m_2^{\tilde {\delta }}B \big ]+ (m_2^{\tilde {\delta }}A m_3-m_3^{\tilde {\delta }}A m_2-Am_2 m_3+Am_3 m_2) \big \} x_1+\ \]
\[ (B-A)^{-1} \big \{ (m_3^{\tilde {\delta }}A-Am_3)(B-A)^{-1} \big [Am_1-m_1^{\tilde {\delta }}B \big ]+ \]
\[ (Am_1-m_1^{\tilde {\delta }}A)(B-A)^{-1}[Am_3-m_3^{\tilde {\delta }}B \big ]+ (m_3^{\tilde {\delta }}A m_1-m_1^{\tilde {\delta }}A m_3-Am_3 m_1+Am_1 m_3) \big \} x_2+ \]
\[ (B-A)^{-1} \big \{ (m_1^{\tilde {\delta }}A-Am_1)(B-A)^{-1} \big [Am_2-m_2^{\tilde {\delta }}B \big ]+ \]
\[ (Am_2-m_2^{\tilde {\delta }}A )(B-A)^{-1}[Am_1-m_1^{\tilde {\delta }}B \big ]+ (m_1^{\tilde {\delta }}A m_2-m_2^{\tilde {\delta }}A m_1-Am_1 m_2+Am_2 m_1) \big \} x_3, 0 \big ). \]

If the loop $L_{\Xi }$ is a group then the Akivis algebra 
${\mathfrak a}_{L_{\Xi }}$ is a Lie algebra.    
The derivation of the condition (II) in Theorem 1 yields 
$m^{\tilde  {\delta }} B=B m$ for all $m \in \mathfrak{m}$.

\smallskip
\noindent
Putting  this in (8)  we obtain for the multiplication in the Lie algebra ${\mathfrak a}_{L_{\Xi }}$ the rule 
\[ [(x_1,m_1),(x_2,m_2)]=( m_1 x_2-m_2 x_1, [m_1,m_2]). \]
The mapping $\gamma:(x,m) \mapsto m: \mathbb R^n \times \mathfrak{m} \to 
\mathfrak{m}$ is an endomorphism from the Akivis algebra 
$\mathfrak{a}_{L_{\Xi }}$ onto 
the Lie algebra $\mathfrak{m}$ since 
\[ \Big [ (x_1,m_1), (x_2,m_2) \Big ] ^{\gamma }= [m_1,m_2] =  \Big [ (x_1,m_1) ^{\gamma }, (x_2,m_2) ^{\gamma } \Big ] \]
and in $\mathfrak{m}$ the Jacobi identity holds. Hence we have the following 

\begin{Prop} The Akivis algebra $\mathfrak{a}_{L_{\Xi }}$ of the loop 
$L_{\Xi }=L_{A,B,\Gamma _0, \delta }$ is a semidirect product $\mathbb R^n 
 \rtimes \mathfrak{m}$ of the commutative Lie algebra 
$\mathbb R^n$ by the Lie algebra  $\mathfrak{m}$ of the group $\Gamma_0$.
\end{Prop}

Let $L_{a,b, \mathbb R^*, \delta}$ be a proper loop constructed in 
$\bf 3.9$ b) of Section 3.  
Then the map $\tilde {\delta }:\mathbb R \to \mathbb R$ corresponding to 
$\delta $ is the automorphism  $x \mapsto (2k+1)x$ with $k \neq 0$. 
We have 
\[ [(x_1,m_1),(x_2,m_2)]= \big ([(2k+1)b -a](b-a)^{-1} (m_1 x_2-m_2 x_1) ,0 \big ) \]
and  
\[ \langle (x_1,m_1), (x_2,m_2), (x_3,m_3) \rangle = \big (-4 k^2 b a (b-a)^{-2} m_1 m_3 x_2 , 0 \big ). \]
Using these expressions  we see that both sides $(9)$
of the Akivis identity  are  equal to $(0,0)$. 

{\it These examples show that there are proper non-connected  smooth loops 
$L_{\Xi }$ of positive dimension having Lie algebras as their Akivis 
algebras.}

\bigskip
\noindent
Mathematisches Institut 
\newline
der Universit\"at Erlangen-N\"urnberg
\newline
Bismarckstr. 1 $\frac{1}{2}$ 
\newline
D-91054 Erlangen
\newline
Germany
\newline
E-MAIL:figula@mi.uni-erlangen.de
\newline
and
\newline
Institute of Mathematics
\newline
University of Debrecen
\newline
P.O.B. 12
\newline
H-4010 Debrecen
\newline
Hungary
\newline 
E-MAIL:figula@math.klte.hu

\bigskip
\noindent
Mathematisches Institut 
\newline
der Universit\"at Erlangen-N\"urnberg
\newline
Bismarckstr. 1 $\frac{1}{2}$ 
\newline
D-91054 Erlangen
\newline
Germany
\newline
E-MAIL:strambach@mi.uni-erlangen.de

\end{document}